\documentclass[11pt]{amsart}
\usepackage{amssymb}
\begin{document}
\theoremstyle{plain}
\newtheorem{thm}{Theorem}[section]
\newtheorem{theorem}[thm]{Theorem}
\newtheorem{lemma}[thm]{Lemma}
\newtheorem{corollary}[thm]{Corollary}
\newtheorem{proposition}[thm]{Proposition}
\newtheorem{conjecture}[thm]{Conjecture}
%%%%%%%%%%%%%%%%%%%% Text roman %%%%%%%%%%%%%%%%%%%%%%%%%%%%%
\theoremstyle{definition}
\newtheorem{construction}[thm]{Construction}
\newtheorem{notations}[thm]{Notations}
\newtheorem{question}[thm]{Question}
\newtheorem{problem}[thm]{Problem}
\newtheorem{remark}[thm]{Remark}
\newtheorem{remarks}[thm]{Remarks}
\newtheorem{definition}[thm]{Definition}
\newtheorem{claim}[thm]{Claim}
\newtheorem{assumption}[thm]{Assumption}
\newtheorem{assumptions}[thm]{Assumptions}
\newtheorem{properties}[thm]{Properties}
\newtheorem{example}[thm]{Example}
\newtheorem{comments}[thm]{Comments}
\newtheorem{blank}[thm]{}
\newtheorem{observation}[thm]{Observation}
\newtheorem{defn-thm}[thm]{Definition-Theorem}

\newcommand{\sM}{{\mathcal M}}

%%%%%%%%%%%%%%%%%%%%%%%%%%%%%%%%%%%%%%%%%%%%%%%%%%%%%%%%%%%%%%

\title[New properties of the intersection numbers]{New properties of the intersection numbers\\ on moduli spaces of curves}
        \author{Kefeng Liu}
        %\dedicatory{Center of Math Sciences, Zhejiang University;
         %   Department of Mathematics, UCLA}
        \address{Center of Mathematical Sciences, Zhejiang University, Hangzhou, Zhejiang 310027, China;
                Department of Mathematics,University of California at Los Angeles,
                Los Angeles, CA 90095-1555, USA}
        \email{liu@math.ucla.edu, liu@cms.zju.edu.cn}
        \author{Hao Xu}
        \address{Center of Mathematical Sciences, Zhejiang University, Hangzhou, Zhejiang 310027, China}
        \email{haoxu@cms.zju.edu.cn}

        \begin{abstract}
            We present certain new properties about the intersection numbers on moduli spaces
            of curves $\overline{\sM}_{g,n}$, including a simple explicit formula of $n$-point functions and several new identities of intersection numbers. In particular
            we prove a new identity, which together with a conjectural
            identity implies the famous Faber's
            conjecture about relations in $\mathcal
            R^{g-2}(\sM_g)$.
            These new identities clarify the mysterious constant in Faber's
            conjecture and uncover novel combinatorial structures of intersection numbers.

            We also discuss some numerical properties of
            Hodge integrals which have provided numerous inspirations for this work.
        \end{abstract}
    \maketitle

\maketitle

\section{Introduction}

Denote by $\overline{\sM}_{g,n}$ the moduli space of stable
$n$-pointed genus $g$ complex algebraic curves. We have the
forgetting the last marked point morphism
$$
\pi: \overline{\sM}_{g,n+1}\longrightarrow \overline{\sM}_{g,n}.
$$
Denote by $\sigma_1,\dots,\sigma_n$ the canonical sections of $\pi$,
and by $D_1,\dots,D_n$ the corresponding divisors in
$\overline{\sM}_{g,n+1}$. Let $\omega_{\pi}$ be the relative
dualizing sheaf, we have the following tautological classes on
moduli spaces of curves.
\begin{align*}
\psi_i&=c_1(\sigma_i^*(\omega_{\pi}))\\
\kappa_i&=\pi_*\left(c_1\left(\omega_{\pi}\left(\sum D_i\right)\right)^{i+1}\right)\\
\lambda_l&=c_l(\pi_*(\omega_{\pi})),\quad 1\leq l\leq g.
\end{align*}
The classes $\kappa_i$ were first introduced by Mumford \cite{Mu} on
$\overline{\sM}_g$, their generalization to $\overline{\sM}_{g,n}$
here is due to Arbarello and Cornalba \cite{Ar-Co}.

We use Witten's notation in this paper,
$$\langle\tau_{d_1}\cdots\tau_{d_n}\lambda_{b_1}\cdots\lambda_{b_k}\rangle_g:=\int_{\overline{\sM}_{g,n}}\psi_1^{d_1}\cdots\psi_n^{d_n}\lambda_{b_1}\cdots\lambda_{b_k}.$$

The moduli space of curves is a central object of study in algebraic
geometry. The intersection theory of tautological classes on the
moduli space of curves is a very important subject and has close
connections to string theory, quantum gravity and many branches of
mathematics.

Intersection numbers involving only $\psi$ classes can be computed
recursively by the celebrated Witten's conjecture \cite{Wi} (proved
by Kontsevich \cite{Ko}) or by the formula of $n$-point functions
\cite{LX}. General intersections involving $\psi, \lambda, \kappa$
or boundary divisor classes can be reduced to intersections of
$\psi$ classes by Faber's algorithm \cite{Fabal}.

The following form of Witten's conjecture is called the DVV formula
\cite{Dij}.
\begin{multline*}
\langle\tau_{k+1}\tau_{d_1}\cdots\tau_{d_n}\rangle_g=\frac{1}{(2k+3)!!}\left[\sum_{j=1}^n
\frac{(2k+2d_j+1)!!}{(2d_j-1)!!}\langle\tau_{d_1}\cdots
\tau_{d_{j}+k}\cdots\tau_{d_n}\rangle_g\right.\\
+\frac{1}{2}\sum_{r+s=k-1}
(2r+1)!!(2s+1)!!\langle\tau_r\tau_s\tau_{d_1}\cdots\tau_{d_n}\rangle_{g-1}\\
\left.+\frac{1}{2}\sum_{r+s=k-1} (2r+1)!!(2s+1)!!
\sum_{\underline{n}=I\coprod J}\langle\tau_r\prod_{i\in
I}\tau_{d_i}\rangle_{g'}\langle\tau_s\prod_{i\in
J}\tau_{d_i}\rangle_{g-g'}\right]
\end{multline*}
where $\underline{n}=\{1,2,\ldots,n\}$.

In 1993, Carel Faber \cite{Fab} proposed his remarkable
conjectures about the structure of tautological ring $\mathcal
R^*(\sM_g)$. In the past decade, Faber's conjecture motivated a
tremendous progress toward understanding of the topology of moduli
spaces of curves. For background materials, we recommend Vakil's
excellent survey \cite{Vak}.

An important part of Faber's conjectures is the famous Faber's
intersection number conjecture, which is the following relations
in $\mathcal R^{g-2}(\sM_g)$, if $\sum_{j=1}^n d_j=g-2$,
\begin{equation}
\sum_{\sigma\in
S_n}\kappa_\sigma=\frac{(2g-3+n)!(2g-1)!!}{(2g-1)!\prod_{j=1}^{n}(2d_j+1)!!}\kappa_{g-2},
\end{equation}
where $\sum_{\sigma\in
S_n}\kappa_\sigma=(\pi_{1}\dots\pi_n)_*(\psi_1^{d_1+1}\dots\psi_n^{d_n+1})$
and $\kappa_\sigma$ is defined as follows. Write the permutation
$\sigma$ as a product of $\nu(\sigma)$ disjoint cycles, including
1-cycles: $\sigma=\beta_1\cdots\beta_{\nu(\sigma)}$, where we
think of the symmetric group $S_n$ as acting on the $n$-tuple
$(d_1,\dots ,d_n)$. Denote by $|\beta|$ the sum of the elements of
a cycle $\beta$. Then $
\kappa_\sigma=\kappa_{|\beta_1|}\kappa_{|\beta_2|}\dots
\kappa_{|\beta_{\nu(\sigma)}|}$.

By the work of Looijenga \cite{Lo} and Faber \cite{Fa1}, we know
that $\mathcal R^{g-2}(\sM_g)=\mathbb Q$ is $1$-dimensional and
$\lambda_g\lambda_{g-1}$ vanishes on the boundary of
$\overline{\sM}_g$. So Faber's intersection number conjecture is
equivalent to the following Hodge integral identity.

If $d_j\geq1$ and $\sum_{j=1}^n(d_j-1)=g-2$,
\begin{eqnarray}
\int_{\overline{\sM}_{g,n}}\psi_1^{d_1}\dots\psi_n^{d_n}\lambda_g\lambda_{g-1}&=&\int_{\overline{\sM}_g}\sum_{\sigma\in\nonumber
S_n}\kappa_\sigma\lambda_g\lambda_{g-1}\\
&=&\frac{(2g-3+n)!(2g-1)!!}{(2g-1)!\prod_{j=1}^{n}(2d_j-1)!!}\int_{\overline{\sM}_g}\kappa_{g-2}\lambda_g\lambda_{g-1}\\
&=&\frac{(2g-3+n)!|B_{2g}|}{2^{2g-1}(2g)!\prod_{j=1}^{n}(2d_j-1)!!}.\nonumber
\end{eqnarray}
where
$\int_{\overline{\sM}_g}\kappa_{g-2}\lambda_g\lambda_{g-1}=\frac{|B_{2g}|(g-1)!}{2^g(2g)!}$
is proved by Faber \cite{Fa1}.

Since $\lambda_{g}\lambda_{g-1}=(-1)^{g-1}(2g-1)!\cdot{\rm
ch}_{2g-1}(\mathbb E)$, we use Mumford's formula \cite{Mu} for the
Chern character of Hodge bundles
$$
{\rm ch}_{2g-1}(\mathbb E)= \frac{ B_{2g}}{(2g)!}\left[
\kappa_{2g-1}-\sum_{i=1}^{n}\psi_i^{2g-1}+\frac12\sum_{\xi\in\Delta}{l_{\xi}}_*
\left(\sum_{i=0}^{2g-2}\psi_{n+1}^i(-\psi_{n+2})^{2g-2-i}\right)\right]
$$
to get
\begin{multline}
\frac{(2g-3+n)!}{2^{2g-1}(2g-1)!}\cdot\frac{1}{\prod_{j=1}^n(2d_j-1)!!}\\
=\langle\tau_{d_1}\cdots\tau_{d_n}\tau_{2g}\rangle_g-\sum_{j=1}^{n}
\langle\tau_{d_1}\cdots\tau_{d_{j-1}}\tau_{d_j+2g-1}\tau_{d_{j+1}}\cdots\tau_{d_n}\rangle_g\\
+\frac{1}{2}\sum_{j=0}^{2g-2}(-1)^j\langle\tau_{2g-2-j}\tau_j\tau_{d_1}\cdots\tau_{d_n}\rangle_{g-1}\\
+\frac{1}{2}\sum_{\underline{n}=I\coprod
J}\sum_{j=0}^{2g-2}(-1)^j\langle\tau_{j}\prod_{i\in
I}\tau_{d_i}\rangle_{g'}\langle\tau_{2g-2-j}\prod_{i\in
J}\tau_{d_i}\rangle_{g-g'}
\end{multline}
where $d_j\geq1$, $\sum_{j=1}^{n}(d_j-1)=g-2$ and
$\underline{n}=\{1,2,\ldots,n\}$.

In fact, the proportional constant in identity (1) is observed
experimentally by Faber from identity (3), through implementing
Witten's conjecture to calculate intersection numbers.

In this paper, we will abuse terminology by calling either of the
equivalent identities (1), (2) and (3) the Faber's conjecture.

Getzler and Pandharipande \cite{Ge-Pa} derive Faber's conjecture
from the degree 0 Virasoro conjecture for $\mathbb P^2$. On the
other hand, Givental \cite{Giv} has announced a proof of Virasoro
conjecture for $\mathbb P^n$. Y.P. Lee and R. Pandharipande are
writing a book supplying the details. Recently, Goulden, Jackson and
Vakil \cite{GJV} have given a more direct and enlightening proof of
Faber's conjecture for up to three points and explained their
approach for the general case. Their method of proof is a marvelous
synthesis of geometry and combinatorics, which has already found
other elegant applications \cite{CLL, GJV2}.

In this paper, we present a series of simple new identities of
intersection numbers, aiming to clarify combinatorial structures in
Faber's conjecture.

The explicit formula of $n$-point functions (as stated in theorem
2.8 and proved in \cite{LX}) for intersection numbers will play an
important role in this work.

\section*{Acknowledgements}
The authors would like to express our special thanks to Ravi Vakil
for many helpful comments and discussions. We wish to thank
Professor Enrico Arbarello, Sergei Lando and Edward Witten for
useful comments and their interests in this work. We thank Professor
Carel Faber for his wonderful Maple program for computing Hodge
integrals and for communicating Zagier's three-point function to us.
We also thank the anonymous referee for detailed comments that have
greatly improved the manuscript.

\section{Several new identities of intersection numbers}

Now we announce the following identity of intersection numbers
which clarifies the mysterious constants in Faber's conjecture.
\begin{theorem}
Let $d_j\geq1$ and $\sum_{j=1}^{n}(d_j-1)=g-1$. Then
\begin{equation}
\sum_{j=0}^{2g}(-1)^j\langle\tau_{2g-j}\tau_j\tau_{d_1}\cdots\tau_{d_n}\rangle_{g}=\frac{(2g-1+n)!}{2^{2g}(2g+1)!}\cdot\frac{1}{\prod_{j=1}^n(2d_j-1)!!}
\end{equation}
\end{theorem}

When $n=1$, identity (4) becomes
$$\sum_{j=0}^{2g}(-1)^j\langle\tau_{2g-j}\tau_j\tau_{g}\rangle_{g}=\frac{1}{2^{2g}(2g+1)!!}$$
which has been proved in \cite{Fa1, Fa-Pa2}.

From theorem 2.1, we see that Faber's conjecture (3) is equivalent
to the following simpler identity.
\begin{conjecture}
Let $d_j\geq0$, $\sum_{j=1}^{n}d_j=g+n-2$ and
$\underline{n}=\{1,2,\ldots,n\}$. Then
\begin{eqnarray}
\langle\tau_{d_1}\cdots\tau_{d_n}\tau_{2g}\rangle_g&=&\sum_{j=1}^{n}
\langle\tau_{d_1}\cdots\tau_{d_{j-1}}\tau_{d_j+2g-1}\tau_{d_{j+1}}\cdots\tau_{d_n}\rangle_g\nonumber\\
&&-\frac{1}{2}\sum_{\underline{n}=I\coprod
J}\sum_{j=0}^{2g-2}(-1)^j\langle\tau_{j}\prod_{i\in
I}\tau_{d_i}\rangle_{g'}\langle\tau_{2g-2-j}\prod_{i\in
J}\tau_{d_i}\rangle_{g-g'}
\end{eqnarray}
\end{conjecture}

Accompanying identity (4) of theorem 2.1, we also have the
following vanishing theorem of intersection numbers.
\begin{theorem}
Let $K>g$, $d_j\geq0$ and $\sum_{j=1}^{n}d_j=3g+n-2K-1$. Then
\begin{equation}
\sum_{j=0}^{2K}(-1)^j\langle\tau_{2K-j}\tau_j\tau_{d_1}\cdots\tau_{d_n}\rangle_g=0
\end{equation}
\end{theorem}

The following corollary of Theorem 2.3 is a complement to Conjecture
2.2.
\begin{corollary}
Let $K>g$, $d_j\geq0$ and $\sum_{j=1}^{n}d_j=3g+n-2K-2$. Then
\begin{eqnarray}
\langle\tau_{d_1}\cdots\tau_{d_n}\tau_{2K}\rangle_g&=&\sum_{j=1}^{n}
\langle\tau_{d_1}\cdots\tau_{d_{j-1}}\tau_{d_j+2K-1}\tau_{d_{j+1}}\cdots\tau_{d_n}\rangle_g\nonumber\\
&&-\frac{1}{2}\sum_{\underline{n}=I\coprod
J}\sum_{j=0}^{2K-2}(-1)^j\langle\tau_{j}\prod_{i\in
I}\tau_{d_i}\rangle_{g'}\langle\tau_{2K-2-j}\prod_{i\in
J}\tau_{d_i}\rangle_{g-g'}
\end{eqnarray}
\end{corollary}
\begin{proof} The identity follows from Mumford's formula \cite{Mu}
$${\rm ch}(\mathbb E)=g+\sum_{k=1}^\infty\frac{B_{2k}}{(2k)!}\left[\kappa_{2k-1}-\sum_{i=1}^n\psi_i^{2k-1}+\frac12\sum_{\xi\in\Delta}{l_{\xi}}_*
\left(\sum_{i=0}^{2k-2}\psi_{n+1}^i(-\psi_{n+2})^{2k-2-i}\right)\right]$$
where ${\rm ch}_{2k-1}(\mathbb E)=0$, when $k>g$.
\end{proof}

Amazingly we also found the following conjectural identity
experimentally. Please compare with the identities (3) and (5),
\begin{conjecture}
Let $g\geq2$, $d_j\geq1$ and $\sum_{j=1}^{n}(d_j-1)=g$. Then
\begin{multline}
\frac{(2g-3+n)!}{2^{2g+1}(2g-3)!}\cdot\frac{1}{\prod_{j=1}^n(2d_j-1)!!}\\
=\langle\tau_{d_1}\cdots\tau_{d_n}\tau_{2g-2}\rangle_g-\sum_{j=1}^{n}
\langle\tau_{d_1}\cdots\tau_{d_{j-1}}\tau_{d_j+2g-3}\tau_{d_{j+1}}\cdots\tau_{d_n}\rangle_g\\
+\frac{1}{2}\sum_{\underline{n}=I\coprod
J}\sum_{j=0}^{2g-4}(-1)^j\langle\tau_{j}\prod_{i\in
I}\tau_{d_i}\rangle_{g'}\langle\tau_{2g-4-j}\prod_{i\in
J}\tau_{d_i}\rangle_{g-g'}.
\end{multline}
\end{conjecture}

Since $(2g-3)!\cdot{\rm ch}_{2g-3}(\mathbb
E)=(-1)^{g-1}(3\lambda_{g-3}\lambda_g-\lambda_{g-1}\lambda_{g-2})$,
it's easy to see that the above identity (8) is equivalent to the
following identity of Hodge integrals,
\begin{conjecture}
Let $g\geq2$, $d_j\geq1$ and $\sum_{j=1}^{n}(d_j-1)=g$. Then
\begin{align}
\frac{2g-2}{|B_{2g-2}|}&\left(\int_{\overline{\sM}_{g,n}}\psi_1^{d_1}\cdots\psi_n^{d_n}\lambda_{g-1}\lambda_{g-2}-3\int_{\overline{\sM}_{g,n}}\psi_1^{d_1}\cdots\psi_n^{d_n}\lambda_{g-3}\lambda_{g}\right)\nonumber\\
&=\frac{1}{2}\sum_{j=0}^{2g-4}(-1)^j\langle\tau_{2g-4-j}\tau_j\tau_{d_1}\cdots\tau_{d_n}\rangle_{g-1}+\frac{(2g-3+n)!}{2^{2g+1}(2g-3)!}\cdot\frac{1}{\prod_{j=1}^n(2d_j-1)!!}
\end{align}
\end{conjecture}

Note that Faber's identity (3) and all of the above identities
(4)-(9) are compatible with the string and dilaton equations, so
$d_j\geq2$ may be assumed when proving these identities. We have
checked the identities of conjecture 2.2 and conjecture 2.5 for all
$g\leq 20$ by computer.

Now we give a proof of Conjecture 2.5 for $n=1$.
$$\frac{g-1}{2^{2g}(2g+1)!!}=\langle\tau_{g+1}\tau_{2g-2}\rangle_g-\langle\tau_{3g-2}\rangle_g+\sum_{j=0}^{2g-4}(-1)^j\langle\tau_j\tau_{g+1}\rangle\langle\tau_{2g-4-j}\rangle.$$

In the last sum of the right hand side, replace $j$ by $3h-g-2$.
We need to prove
$$\sum_{h=1}^g\frac{(-1)^{g-h}}{24^{g-h}(g-h)!}\langle\tau_{3h-g-2}\tau_{g+1}\rangle_h=\frac{g-1}{2^{2g}(2g+1)!!}+\langle\tau_{3g-2}\rangle_g.$$

Apply the string equation twice, we have
\begin{multline*}
\sum_{h=1}^g\frac{(-1)^{g-h}}{24^{g-h}(g-h)!}\langle\tau_{3h-g-2}\tau_{g+1}\rangle_h\\
=\sum_{h=1}^g\frac{(-1)^{g-h}}{24^{g-h}(g-h)!}\left(\langle\tau_0\tau_{3h-g-1}\tau_{g+1}\rangle_h-\langle\tau_{3h-g-1}\tau_{g}\rangle_h\right)\\
=\sum_{h=1}^g\frac{(-1)^{g-h}}{24^{g-h}(g-h)!}\left(\langle\tau_0\tau_{3h-g-1}\tau_{g+1}\rangle_h-\langle\tau_0\tau_{3h-g}\tau_{g}\rangle_h+\langle\tau_{3h-g}\tau_{g-1}\rangle_h\right).
\end{multline*}

Here we need some knowledge of $n$-point functions such as in the
papers \cite{Fa1, Fa-Pa2, LX}. Since
$\sum_{h=0}^g\frac{(-1)^{g-h}}{24^{g-h}(g-h)!}\langle\tau_0\tau_{3h-g-1}\tau_{g+1}\rangle_h$
and
$\sum_{h=0}^g\frac{(-1)^{g-h}}{24^{g-h}(g-h)!}\langle\tau_0\tau_{3h-g}\tau_{g}\rangle_h$
are respectively the coefficient of $y^{2g-1}z^{g+1}$ and
$y^{2g}z^{g}$ in
$$\exp\left(\frac{z^3}{24}\right)\sum_{k\geq0}\frac{k!}{(2k+1)!}\left(\frac{1}{2}yz(y+z)\right)^k,$$
we have
\begin{align*}
\sum_{h=1}^g\frac{(-1)^{g-h}}{24^{g-h}(g-h)!}\langle\tau_0\tau_{3h-g-1}\tau_{g+1}\rangle_h&=\frac{g!}{(2g+1)!}\cdot\frac{g}{2^g},\\
\sum_{h=1}^g\frac{(-1)^{g-h}}{24^{g-h}(g-h)!}\langle\tau_0\tau_{3h-g}\tau_{g}\rangle_h&=\frac{g!}{(2g+1)!}\cdot\frac{1}{2^g}.
\end{align*}

Moreover, it has been proved in \cite{Fa1, Fa-Pa2} that
$$\sum_{h=1}^g\frac{(-1)^{g-h}}{24^{g-h}(g-h)!}\langle\tau_{3h-g}\tau_{g-1}\rangle_h=\frac{1}{24^gg!},$$
so we conclude the proof.

\bigskip
Both of Theorem 2.1 and Theorem 2.3 follow from our recently
obtained $n$-point functions for intersection numbers.

\begin{definition}
We call the following generating function
$$F(x_1,\dots,x_n)=\sum_{g=0}^{\infty}\sum_{\sum d_j=3g-3+n}\langle\tau_{d_1}\cdots\tau_{d_n}\rangle_g\prod_{j=1}^n x_j^{d_j}$$
the $n$-point function.
\end{definition}

Note that the left hand side of identity (4) in Theorem 2.1 is
$$\left[F(y,-y,x_1,\dots,x_n)\right]_{y^{2g-2}\prod_{j=1}^{n} x_j^{d_j}},$$
which is the coefficient of the monomial $y^{2g-2}\prod_{j=1}^{n}
x_j^{d_j}$ in the special $(n+2)$-point function
$F(y,-y,x_1,\dots,x_n)$.

It's not an easy task to get explicit formulae for $n$-point
functions \cite{BH}. Okounkov \cite{Ok} has obtained a marvelous
analytic formula for $n$-point functions, however it seems very
difficult to extract information of coefficients from this analytic
formula. What's more interesting is to find some well organized
series expansion for the $n$-point function.

We introduce the following ``normalized'' $n$-point function
$$G(x_1,\dots,x_n)=\exp\left(\frac{-\sum_{j=1}^n
x_j^3}{24}\right)\cdot F(x_1,\dots,x_n).$$

\begin{theorem} \cite{LX} For $n\geq2$,
\begin{equation*}
G(x_1,\dots,x_n)=\sum_{r,s\geq0}\frac{(2r+n-3)!!}{4^s(2r+2s+n-1)!!}P_r(x_1,\dots,x_n)\Delta(x_1,\dots,x_n)^s,
\end{equation*}
where $P_r$ and $\Delta$ are homogeneous symmetric polynomials
defined by
\begin{align*}
\Delta(x_1,\dots,x_n)&=\frac{(\sum_{j=1}^nx_j)^3-\sum_{j=1}^nx_j^3}{3},\\
P_r(x_1,\dots,x_n)&=\left(\frac{1}{2\sum_{j=1}^nx_j}\sum_{\underline{n}=I\coprod
J}(\sum_{i\in I}x_i)^2(\sum_{i\in J}x_i)^2 G(x_I)
G(x_J)\right)_{3r+n-3}\\
&=\frac{1}{2\sum_{j=1}^nx_j}\sum_{\underline{n}=I\coprod
J}(\sum_{i\in I}x_i)^2(\sum_{i\in J}x_i)^2\sum_{r'=0}^r
G_{r'}(x_I)G_{r-r'}(x_J),
\end{align*}
where $I,J\ne\emptyset$, $\underline{n}=\{1,2,\ldots,n\}$ and
$G_g(x_I)$ denotes the degree $3g+|I|-3$ homogeneous component of
the normalized $|I|$-point function $G(x_{k_1},\dots,x_{k_{|I|}})$,
where $k_j\in I$.
\end{theorem}

The above formula generalizes Dijkgraaf's two-point function
\cite{Fa1} and Zagier's three-point function \cite{Za} obtained more
than ten years ago.

Let's turn to the normalized special $(n+2)$-point function,
$$G(y,-y,x_1,\dots,x_n)=\exp\left(\frac{-\sum_{j=1}^n
x_j^3}{24}\right)\cdot F(y,-y,x_1,\dots,x_n).$$

We have the following theorem \cite{LX} about the coefficients of
$G(y,-y,x_1,\dots,x_n)$, which is just a reformulation of Theorem
2.3 and Theorem 2.1.

\begin{theorem}
Let $g\geq0$ and $n\geq1$. We have
\begin{enumerate}
\item Let $K>g, d_j\geq0$ and $\sum_{j=1}^n d_j=3g-1+n-2K$. Then
$$\left[G(y,-y,x_1,\dots,x_n)\right]_{y^{2K}\prod_{j=1}^n x_j^{d_j}}=0.$$

\item Let $d_j\geq1$ and $\sum_{j=1}^n d_j=g-1+n$. Then
$$\left[G(y,-y,x_1,\dots,x_n)\right]_{y^{2g}\prod_{j=1}^n x_j^{d_j}}=\frac{(2g+n-1)!}{4^{g}(2g+1)!\cdot\prod_{j=1}^n(2d_j-1)!!}.$$
\end{enumerate}
\end{theorem}

\section{Vanishing phenomenon of intersection numbers}
We investigate a sort of vanishing phenomenon of intersection
numbers and seek further clarification for our simplified version of
Faber's conjecture. First we state a generalization of Theorem 2.3.
\begin{conjecture}
Let $K>g$, $d_j\geq0$ and $\Lambda$ be a monomial of $\lambda$
classes. Then
\begin{equation*}
\sum_{j=0}^{2K}(-1)^j\langle\tau_{2K-j}\tau_j\Lambda\prod_{j=1}^n\tau_{d_j}\rangle_g=0
\end{equation*}
\end{conjecture}

Take $\Lambda=\lambda_1$ in Conjecture 3.1. Since $\lambda_1={\rm
ch}_1(\mathbb E)$, we have
\begin{multline*}
0=12
\sum_{j=0}^{2K}(-1)^j\langle\tau_{2K-j}\tau_j\lambda_1\prod_{j=1}^n\tau_{d_j}\rangle_g\\
=-\sum_{j=0}^{2K}(-1)^j\left(\langle\tau_{2K-j+1}\tau_j\prod_{j=1}^n\tau_{d_j}\rangle_g+\langle\tau_{2K-j}\tau_{j+1}\prod_{j=1}^n\tau_{d_j}\rangle_g\right)\\
+\sum_{\underline{n}=I\coprod
J}\sum_{j=0}^{2K}(-1)^j\langle\tau_j\tau_0\prod_{i\in
I}\tau_{d_i}\rangle_{g'}\langle\tau_{2K-j}\tau_0\prod_{i\in
J}\tau_{d_i}\rangle_{g-g'}\\
=\sum_{\underline{n}=I\coprod
J}\sum_{j=0}^{2K}(-1)^j\langle\tau_j\tau_0\prod_{i\in
I}\tau_{d_i}\rangle_{g'}\langle\tau_{2K-j}\tau_0\prod_{i\in
J}\tau_{d_i}\rangle_{g-g'}-2\langle\tau_{2K+1}\tau_0\prod_{j=1}^n\tau_{d_j}\rangle_g.
\end{multline*}
We have used Theorem 2.3 in the above equations. In fact, we have
identities much more general than the above equation.

\begin{conjecture} We have
\begin{enumerate}
\item
Let $K\geq g$, $r,s\geq0$, $d_j\geq0$ and
$\sum_{j=1}^{n}d_j=3g+n-2K-r-s-2$. Then
\begin{multline*}
\langle\tau_{2K+r+1}\tau_s\prod_{j=1}^n\tau_{d_j}\rangle_g+\langle\tau_{2K+s+1}\tau_r\prod_{j=1}^n\tau_{d_j}\rangle_g\\
=\sum_{\underline{n}=I\coprod
J}\sum_{j=0}^{2K}(-1)^j\langle\tau_j\tau_r\prod_{i\in
I}\tau_{d_i}\rangle_{g'}\langle\tau_{2K-j}\tau_s\prod_{i\in
J}\tau_{d_i}\rangle_{g-g'}.
\end{multline*}

\item
Let $r,s\geq0$, $d_j\geq1$ and $\sum_{j=1}^{n}d_j=g+n-r-s$. Then
\begin{multline*}
\frac{1}{(2r+1)!!(2s+1)!!}\cdot\frac{(2g-1+n)!}{4^g(2g-1)!\prod_{j=1}^n(2d_j-1)!!}\\
=\langle\tau_{2g+r-1}\tau_s\prod_{j=1}^n\tau_{d_j}\rangle_g+\langle\tau_{2g+s-1}\tau_r\prod_{j=1}^n\tau_{d_j}\rangle_g\\
-\sum_{\underline{n}=I\coprod
J}\sum_{j=0}^{2g-2}(-1)^j\langle\tau_j\tau_r\prod_{i\in
I}\tau_{d_i}\rangle_{g'}\langle\tau_{2g-2-j}\tau_s\prod_{i\in
J}\tau_{d_i}\rangle_{g-g'}.
\end{multline*}
\end{enumerate}
\end{conjecture}

Identities in Conjecture 3.2 have the same structures as results in
the last section, we believe there is a uniform way to prove these
conjectural identities.

Taking $\Lambda={\rm ch}_{2r+1}(\mathbb E)$ in Conjecture 3.1 and
$K>g$, we have
\begin{multline*}
0=\frac{(2r+2)!}{B_{2r+2}}
\sum_{j=0}^{2K}(-1)^j\langle\tau_{2K-j}\tau_j{\rm ch}_{2r+1}(\mathbb E)\prod_{i=1}^n\tau_{d_i}\rangle_g\\
=-\sum_{j=0}^{2K}(-1)^j\left(\langle\tau_{2K-j+2r+1}\tau_j\prod_{i=1}^n\tau_{d_i}\rangle_g+\langle\tau_{2K-j}\tau_{j+2r+1}\prod_{i=1}^n\tau_{d_i}\rangle_g\right)\\
+\sum_{\underline{n}=I\coprod
J}\sum_{i=0}^{2r}\sum_{j=0}^{2K}(-1)^{i+j}\langle\tau_i\tau_j\prod_{t\in
I}\tau_{d_t}\rangle_{g'}\langle\tau_{2r-i}\tau_{2K-j}\prod_{t\in
J}\tau_{d_t}\rangle_{g-g'}\\
=-2\sum_{i=0}^{2r}(-1)^i\langle\tau_{2K+2r+1-i}\tau_i\prod_{j=1}^n\tau_{d_j}\rangle_g\\
+\sum_{\underline{n}=I\coprod
J}\sum_{i=0}^{2r}\sum_{j=0}^{2K}(-1)^{i+j}\langle\tau_i\tau_j\prod_{t\in
I}\tau_{d_t}\rangle_{g'}\langle\tau_{2r-i}\tau_{2K-j}\prod_{t\in
J}\tau_{d_t}\rangle_{g-g'}.
\end{multline*}

It's not difficult to see that conjecture 3.2(1) implies Conjecture
3.1 in the case $\Lambda={\rm ch}_{2r+1}(\mathbb E)$.

Now we present a hierarchy of conjectural identities of intersection
numbers, which provide further insights to Faber's conjecture.

\begin{conjecture} For $m\geq 2$, we have
\begin{enumerate}
\item
Let $K\geq g+\lfloor\frac{m}{2}\rfloor-1$, $r_p\geq 0$, $d_j\geq0$
and $\sum_{j=1}^{n}d_j=3g+n-2K-\sum_{p=1}^m r_p+m-4$. Then
\begin{multline*}
\langle\tau_{2K+2}\prod_{p=1}^m\tau_{r_p}\prod_{j=1}^n\tau_{d_j}\rangle_g=\sum_{j=1}^n\langle\tau_{d_1}\dots\tau_{d_{j-1}}\tau_{d_j+2K+1}\tau_{d_{j+1}}\prod_{p=1}^m\tau_{r_p}\rangle_g\\
-\sum_{\underline{n}=I\coprod
J}\sum_{j=0}^{2K}(-1)^j\langle\tau_j\prod_{p=1}^m\tau_{r_p}\prod_{i\in
I}\tau_{d_i}\rangle_{g'}\langle\tau_{2K-j}\prod_{i\in
J}\tau_{d_i}\rangle_{g-g'}.
\end{multline*}

\item
Let $K=g+\lfloor\frac{m}{2}\rfloor-2$, $d_j\geq1$ and
$\sum_{j=1}^{n}d_j=g+n-\sum_p r_p+m-2\lfloor\frac{m}{2}\rfloor$.
Define
$$
C(g,m,r_p)=
\begin{cases} (\sum_p 2r_p+m)(g+\frac{m-3}{2})& \text{if $m$ is odd,}
\\
\ 1 & \text{if $m$ is even.}
\end{cases}
$$
Then we have
\begin{multline*}
\frac{C(g,m,r_p)}{\prod_{p=1}^m(2r_p+1)!!}\cdot\frac{(2g-3+n+m)!}{4^g(2g-3+m)!\prod_{j=1}^n(2d_j-1)!!}\\
=\langle\tau_{2K+2}\prod_{p=1}^m\tau_{r_p}\prod_{j=1}^n\tau_{d_j}\rangle_g-\sum_{j=1}^n\langle\tau_{d_1}\dots\tau_{d_{j-1}}\tau_{d_j+2K+1}\tau_{d_{j+1}}\prod_{p=1}^m\tau_{r_p}\rangle_g\\
+\sum_{\underline{n}=I\coprod
J}\sum_{j=0}^{2K}(-1)^j\langle\tau_j\prod_{p=1}^m\tau_{r_p}\prod_{i\in
I}\tau_{d_i}\rangle_{g'}\langle\tau_{2K-j}\prod_{i\in
J}\tau_{d_i}\rangle_{g-g'}.
\end{multline*}
If $m$ is odd, we require $r_p\geq1$ in the above identity.
\end{enumerate}
\end{conjecture}

\begin{conjecture} For $m\geq 2$, we have
\begin{enumerate}
\item
Let $K\geq g+\lfloor\frac{m-1}{2}\rfloor$, $s\geq0$, $r_p\geq 0$,
$d_j\geq0$ and $\sum_{j=1}^{n}d_j=3g+n-2K-s-\sum_{p=1}^m r_p+m-3$.
Then
\begin{multline*}
\langle\tau_{2K+s+1}\prod_{p=1}^m\tau_{r_p}\prod_{j=1}^n\tau_{d_j}\rangle_g\\
=\sum_{\underline{n}=I\coprod
J}\sum_{j=0}^{2K}(-1)^j\langle\tau_j\prod_{p=1}^m\tau_{r_p}\prod_{i\in
I}\tau_{d_i}\rangle_{g'}\langle\tau_{2K-j}\tau_s\prod_{i\in
J}\tau_{d_i}\rangle_{g-g'}.
\end{multline*}

\item
Let $K=g+\lfloor\frac{m-1}{2}\rfloor-1$, $d_j\geq1$ and
$\sum_{j=1}^{n}d_j=g+n-s-\sum_p
r_p+m-2\lfloor\frac{m-1}{2}\rfloor-1$. Define
$$
C(g,m,s,r_p)=
\begin{cases} (\sum_p 2r_p-2s+m-1)(g+\frac{m}{2}-1)& \text{if $m$ is even,}
\\
\ 1 & \text{if $m$ is odd.}
\end{cases}
$$
Then we have
\begin{multline*}
\frac{C(g,m,s,r_p)}{(2s+1)!!\prod_{p=1}^m(2r_p+1)!!}\cdot\frac{(2g-2+n+m)!}{4^g(2g-2+m)!\prod_{j=1}^n(2d_j-1)!!}\\
=\langle\tau_{2K+s+1}\prod_{p=1}^m\tau_{r_p}\prod_{j=1}^n\tau_{d_j}\rangle_g\\
-\sum_{\underline{n}=I\coprod
J}\sum_{j=0}^{2K}(-1)^j\langle\tau_j\prod_{p=1}^m\tau_{r_p}\prod_{i\in
I}\tau_{d_i}\rangle_{g'}\langle\tau_{2K-j}\tau_s\prod_{i\in
J}\tau_{d_i}\rangle_{g-g'}.
\end{multline*}
If $m$ is even, we require $s\geq1$ and $r_p\geq1$ in the above
identity.
\end{enumerate}
\end{conjecture}

The following conjecture generalizes results of \cite{LX}.
\begin{conjecture} For $m,l\geq 2$, we have
\begin{enumerate}
\item Let $K>2g+m+l-4$, $r_p,s_p\geq0$, $d_j\geq
0$ and $\sum_{j=1}^{n}d_j=3g+n+m+l-K-\sum_{p=1}^m r_p-\sum_{p=1}^l
s_p-4$. Then
$$
\sum_{\underline{n}=I\coprod
J}\sum_{j=0}^{K}(-1)^j\langle\tau_j\prod_{p=1}^m\tau_{r_p}\prod_{t\in
I}\tau_{d_t}\rangle\langle\tau_{K-j}\prod_{p=1}^l\tau_{s_p}\prod_{t\in
J}\tau_{d_t}\rangle=0.
$$

\item Let $K=2g+m+l-4$, $d_j\geq 1$ and $\sum_{j=1}^{n}d_j=g+n-\sum_p r_p-\sum_p s_p$. Then
\begin{multline*}
\sum_{\underline{n}=I\coprod
J}\sum_{j=0}^{K}(-1)^j\langle\tau_j\prod_{p=1}^m\tau_{r_p}\prod_{t\in
I}\tau_{d_t}\rangle\langle\tau_{K-j}\prod_{p=1}^l\tau_{s_p}\prod_{t\in
J}\tau_{d_t}\rangle\\
=\frac{1}{\prod_{p=1}^m(2r_p+1)!!\prod_{p=1}^l(2s_p+1)!!}\cdot\frac{(-1)^m(2g+n+m+l-3)!}{4^g(2g+m+l-3)!\prod_{j=1}^n(2d_j-1)!!}.
\end{multline*}

\end{enumerate}
\end{conjecture}

\section{Denominators of intersection numbers}

Let {\it denom}$(r)$ denotes the denominator of a rational number
$r$ in reduced form (coprime numerator and denominator, positive
denominator). We define
$$D_{g,n}=lcm\left\{denom\left(\int_{\overline{\sM}_{g,n}}\psi_1^{d_1}\cdots\psi_n^{d_n}\right)\Big{|}\
\sum_{i=1}^{n}d_i=3g-3+n\right\}$$ and for $g\geq 2$,
$$\mathcal D_g=lcm\left\{denom\left(\int_{\overline{\sM}_{g}}\kappa_{a_1}\cdots\kappa_{a_m}\right)\Big{|}\ \sum_{i=1}^{m}a_m=3g-3\right\}$$
where {\it lcm} is the abbreviation of {\it least common multiple}.

We proved in our previous paper \cite{Liu-Xu} that
$D_{g,n}\mid\mathcal D_g$ and $D_{g,n}=\mathcal D_g$ for $n\geq
3g-3$.

Now we present the conjectural exact values of $\mathcal D_g$.

\begin{conjecture}
Let $p$ denotes a prime number and $g\geq 2$. Let ${\rm ord}(p,n)$
denotes the maximum integer such that $p^{{\rm ord}(p,n)}\mid n$.
Then
\begin{enumerate}
\item ${\rm ord}(2,\mathcal D_g)=3g+{\rm ord}(2,g!)$,

\item ${\rm ord}(3,\mathcal D_g)=g+{\rm ord}(3,g!)$,

\item ${\rm ord}(p,\mathcal D_g)=\lfloor\frac{2g}{p-1}\rfloor$ for
$p\geq5$, where $\lfloor x\rfloor$ denotes the maximum integer
that is not larger than $x$.
\end{enumerate}

We order all Witten-Kontsevich tau functions of given genus $g$ by
the following lexicographical rule,
$$\langle\tau_{d_1}\cdots\tau_{d_n}\rangle_g\prec\langle\tau_{k_1}\cdots\tau_{k_m}\rangle_{g}$$
if $n<m$ or $n=m$ and there exists some $i$, such that $d_j=k_j\
\text{for}\ j<i$ and $d_i<k_i$.

If $5\leq p\leq 2g+1$ is a prime number, then the smallest tau
function of genus $g$ in the above lexicographical order that
satisfies ${\rm
ord}(p,denom\langle\tau_{d_1}\cdots\tau_{d_n}\rangle_g)=\lfloor\frac{2g}{p-1}\rfloor$
is
$$\langle\underbrace{\tau_{\frac{p-1}{2}}\cdots\tau_{\frac{p-1}{2}}}_{\lfloor\frac{2g}{p-1}\rfloor}\tau_{d}\rangle_g$$
where
$d+\frac{p-1}{2}\lfloor\frac{2g}{p-1}\rfloor=3g-2+\lfloor\frac{2g}{p-1}\rfloor$.
\end{conjecture}
We have checked Conjecture 4.1 for all $g\leq 20$ by a computer.

\begin{corollary}
We have $D_{g,n}=\mathcal D_g$ for
$n\geq\lfloor\frac{g}{2}\rfloor+1$.
\end{corollary}

\begin{corollary}
Let $\mathcal D_0=1$ and $\mathcal D_1=24$, then $\mathcal
D_{g}\mathcal D_{h}\mid\mathcal D_{g+h}$, for $g,h\geq 0$.
\end{corollary}

Note that $\lfloor\frac{g}{2}\rfloor+1$ is just the number of
codimension one boundary strata of $\overline{\sM}_g$, we don't
know whether this has any implications.

We remark that that $\mathcal D_g$ does not control the denominators
of general Hodge integrals, since we have ${\rm
ord}(5,denom\langle\tau_{19}\lambda_9\rangle_{10})=6>{\rm
ord}(5,\mathcal D_{10})$.

Let $\mathcal S_g$ be the least common multiple of
$\{|Aut(\Sigma_g)|\}$, where $\Sigma_g$ takes over all stable curves
of genus $g$. By arranging components of stable curves in a most
symmetric way, it's not difficult to see that
\begin{align*}
&{\rm ord}(2, \mathcal S_g)\geq 2g+\lfloor g/2\rfloor+\lfloor\lfloor
g/2\rfloor/2\rfloor+\lfloor\lfloor\lfloor
g/2\rfloor/2\rfloor/2\rfloor+\cdots\\
&{\rm ord}(p, \mathcal S_g)\geq k+\lfloor k/p\rfloor+\lfloor\lfloor
k/p\rfloor/p\rfloor+\lfloor\lfloor\lfloor
k/p\rfloor/p\rfloor/p\rfloor+\cdots,\ \text{if\ prime}\ p\geq3,
\end{align*}
where $k=\lfloor\frac{2g}{p-1}\rfloor$. We conjecture that the above
relations are actually equalities giving exact values of $\mathcal
S_g$.

It's easy to see that ${\rm ord}(2, \mathcal D_g)>{\rm ord}(2,
\mathcal S_g),\ {\rm ord}(3, \mathcal D_g)\geq{\rm ord}(3, \mathcal
S_g)$ and for $p\geq5$, ${\rm ord}(p, \mathcal D_g)\leq{\rm ord}(p,
\mathcal S_g)$. Roughly speaking, this means that top intersections
of kappa classes fail to detect all singularities on the orbifold
$\overline{\sM}_g$.

\vskip 30pt
\section{Numerical properties of intersection numbers}

From Okounkov's analytic formula of $n$-point functions \cite{Ok},
we have
$$F(x_1,\dots,x_n)=\sum_{g=0}^{\infty}\sum_{\sum d_j=3g-3+n}\langle\tau_{d_1}\cdots\tau_{d_n}\rangle_g\prod_{j=1}^n x_j^{d_j}<\infty.$$
for arbitrary positive real numbers $x_i$. So
$\langle\tau_{d_1}\cdots\tau_{d_n}\rangle_g$ decreases very rapidly
when $g$ increases. In this section, we will discuss a kind of
multinomial-type property for intersection numbers.
\begin{conjecture}
Let $\Lambda$ be a monomial of the form
$\lambda_1^{k_1}\cdots\lambda_g^{k_g}$. Then for
$\sum_{j=1}^{n}d_j=3g-3+n-\sum_{j=1}^{g}k_j\cdot j$ and $d_1<d_2$,
we have
$$\int_{\overline{\sM}_{g,n}}\psi_1^{d_1}\psi_2^{d_2}\cdots\psi_{n}^{d_n}\Lambda\leq\int_{\overline{\sM}_{g,n}}\psi_1^{d_1+1}\psi_2^{d_2-1}\cdots\psi_{n}^{d_n}\Lambda.$$
\end{conjecture}
Namely the more evenly $3g-3+n$ be distributed among indices, the
larger the value of Hodge integrals.

From the argument of Proposition 5.1 of \cite{Liu-Xu}, we see it's
enough to check only those Hodge integrals with $d_3\geq
2,\ldots,d_n\geq 2$. We have checked Conjecture 5.1 in various
cases.

For $\Lambda=1$, namely in the case of tau functions, we have
checked Conjecture 5.1 for $g\leq 20$. Moreover, for $n=2$, we have
checked all $g\leq 1000$ (using Dijkgraaf's $2$-point function); for
$n=3$, we have checked all $g\leq 100$ (using Zagier's $3$-point
function).

For $\Lambda=\lambda_g$, we have the $\lambda_g$ theorem proved by
Faber and Pandharipande \cite{Fa-Pa},
$$\int_{\overline{\sM}_{g,n}}\psi_1^{d_1}\cdots\psi_{n}^{d_n}\lambda_g=\binom{2g+n-3}{d_1,\dots,d_n}\frac{2^{2g-1}-1}{2^{2g-1}}\frac{|B_{2g}|}{(2g)!}.
$$

For $\Lambda=\lambda_{g}\lambda_{g-1}$, we have the Faber's
conjecture.

For $\Lambda=\lambda_1^{k_1}\cdots\lambda_g^{k_g}$, where
$\sum_{i=1}^gik_i=3g-3$, we can use the same argument as Proposition
5.1 in \cite{Liu-Xu}.

For $\Lambda=\lambda_{g-1}$, there is a closed formula for Hodge
integrals with $\lambda_{g-1}$ class in \cite{Ge-OP}, we have
checked the case of two-point Hodge integrals for $g\leq 100$.

Based on a large amount of experiment, we speculate that the
following generalization of Conjecture 5.1 should be true.

\begin{conjecture}
Consider the following general intersection numbers,
$$\langle\tau_{\underline{d}}\kappa_{\underline{a}}\lambda_{\underline{b}}\rangle_{g,n}:=\int_{\overline{\mathcal{M}}_{g,n}}\psi_{1}^{d_{1}}\cdots\psi_{n}^{d_{n}}\kappa_{a_1}\cdots\kappa_{a_m}\lambda_{b_{1}}
    \cdots\lambda_{b_{k}}$$
where the indices \{$d_j, a_j, b_j$\} are nonnegative integers. If
$p<q$ we have
\begin{align*}
\langle\tau_{p}\tau_{q}\tau_{\underline{d}}\kappa_{\underline{a}}\lambda_{\underline{b}}\rangle_{g,n}&\leq\langle\tau_{p+1}\tau_{q-1}\tau_{\underline{d}}\kappa_{\underline{a}}\lambda_{\underline{b}}\rangle_{g,n},\\
\langle\kappa_{p}\kappa_{q}\tau_{\underline{d}}\kappa_{\underline{a}}\lambda_{\underline{b}}\rangle_{g,n}&\leq\langle\kappa_{p+1}\kappa_{q-1}\tau_{\underline{d}}\kappa_{\underline{a}}\lambda_{\underline{b}}\rangle_{g,n},\\
\langle\lambda_{p}\lambda_{q}\tau_{\underline{d}}\kappa_{\underline{a}}\lambda_{\underline{b}}\rangle_{g,n}&\leq\langle\lambda_{p+1}\lambda_{q-1}\tau_{\underline{d}}\kappa_{\underline{a}}\lambda_{\underline{b}}\rangle_{g,n}.
\end{align*}
\end{conjecture}

We make some remarks about Conjecture 5.2. First recall the
definition of Schur polynomials. Let $E$ be a vector bundle of rank
$r$ on a projective variety $X$ and $\mu$ be a partition of integer
$n$ into integers $\leq r$,
$$r\geq\mu_1\geq\mu_2\geq\dots\geq\mu_n\geq0.$$
Define the Schur
polynomials $s_\mu(E)$ by
$$s_\mu(E)=\det(c_{\mu_i+j-i}(E))_{1\leq i,j\leq n}.$$
So $s_\mu(E)$ is a polynomial in Chern classes of $E$ with
weighted degree $n$. If $E$ is ample, then we know that $s_\mu(E)$
behave numerical positively on $X$.

However the Hodge bundle $\mathbb E$ on $\overline{\sM}_g$ is not
ample, only its determinant $\det(\mathbb E)$ is semi-ample. If
$\mu=(b_2,b_1,0,\dots,0)$, then $s_\mu(\mathbb
E)=\lambda_{b_2}\lambda_{b_1}-{\lambda_{b_2+1}\lambda_{b_1-1}}$. So
conjecture 5.2(3) essentially says that Schur polynomials of the
form $s_{(b_2,b_1)}(\mathbb E)$ behave numerical positively on the
moduli space of curves.

There are several natural consequences from Conjecture 5.2. For
example, we could get simple lower and upper bounds for
intersections of kappa classes, the so called higher Weil-Petersson
volumes of the moduli space of curves \cite{KMZ}.

\begin{corollary} Let $a_j\geq0$, $\sum_{j=1}^{m}a_j=3g-3+n$ and
$g\geq1$. We have
$$\frac{(2g-2+n)^{m-1}}{24^g\cdot g!}\leq\langle\kappa_{a_1}\cdots\kappa_{a_m}\rangle_{g,n}\leq\frac{\langle\kappa_1^{3g-3+n}\rangle_{g,n}}{(2g-2+n)^{3g-3+n-m}}.$$
In particular, we get a simple lower bound for Weil-Petersson
volumes
$$\langle\kappa_1^{3g-3+n}\rangle_{g,n}\geq\frac{(2g-2+n)^{3g-4+n}}{24^g\cdot g!}.$$
\end{corollary}

\begin{proof} Since $\kappa_0=2g-2+n$ in $\mathcal R^0(\overline{\sM}_{g,n})$,
we have
$$\langle\kappa_{a_1}\cdots\kappa_{a_m}\rangle_{g,n}\geq\langle\kappa_{3g-3+n}\kappa_{0}^{m-1}\rangle_{g,n}=\frac{(2g-2+n)^{m-1}}{24^g\cdot g!}.$$
The other inequality can be proved similarly.
\end{proof}

\begin{corollary}
For $d_j\geq0$ and $\sum_{j=1}^{n}d_j=3g-3+n$,
$$\langle\tau_{d_1}\cdots\tau_{d_n}\rangle_g\geq\frac{1}{24^g\cdot g!}.$$
\end{corollary}

%$$ \ \ \ \ $$

\end{document}